\documentclass[11pt]{article}

\usepackage{amsmath}
\usepackage{amssymb}

\usepackage{latexsym}
\usepackage[usenames]{color}
\def\@begintheorem#1#2{\par\bgroup{\sc #1\ #2. }\it\ignorespaces}
\def\@opargbegintheorem#1#2#3{\par\bgroup{\sc #1\ #2\ (#3). }\it\ignorespaces}
\def\@endtheorem{\egroup}

\newtheorem{theorem}{Theorem}

\newtheorem{corollary}{Corollary}

\def\@abssec#1{\vspace{.05in}\footnotesize \parindent .2in
{\bf #1. }\ignorespaces}
\def\abstract{\@abssec{Abstract}}
\def\keywords{\@abssec{Key words}}
\def\AMSMOS{\@abssec{AMS(MOS) subject classifications}}

\title{\bf Generalizations of Sylvester's determinantal identity}

\begin{document}

\author{ A. Karapiperi
\\ {\tt akarapi@math.unipd.it}
\and M.
Redivo--Zaglia\thanks{Corresponding author:
Universit\`a di Padova, Dipartimento di Matematica,
       Via Trieste, 63 - 35121 Padova - Italy,
        Email: Michela.RedivoZaglia@unipd.it}
\\  {\tt Michela.RedivoZaglia@unipd.it}
 \and M.R. Russo
\\  {\tt mrrusso@math.unipd.it}
\\ \\
 Dipartimento di Matematica \\ Universit\`a degli Studi di Padova, Italy.  }

\date{~}

\maketitle

\begin{abstract}
In this paper we deal with the noteworthy Sylvester's determinantal
identity and some of its generalizations. We report the
formulae due to Yakovlev,  to Gasca, Lopez--Carmona, Ramirez,
to Beckermann, Gasca, M\"uhlbach, and to Mulders in a
unified formulation which allows to understand them better and to compare them. Then, we propose a
different generalization of Sylvester's classical formula. This new
generalization expresses the determinant of a matrix in relation
with the determinant of the bordered matrices obtained adding more
than one row and one column to the original matrix. Sylvester's identity is recovered as a particular case.

\end{abstract}

\section{Introduction}

Sylvester's determinantal identity \cite{Sylv} is a well-known
identity in matrix analysis which expresses a determinant composed
of bordering determinants in terms of the original one. It has been
extensively studied, both in the algebraic and in the combinatorial
context, and it is usually written in the form used by Bareiss
\cite{Bair}.

In the last years some generalizions of this identity have
been proposed, with useful applications mainly in interpolation and extrapolation
problems \cite{THEBOOK,schur}. One of the first generalizations seems to have been proposed
by Yakovlev \cite{Yak1,Yak2} who used it in the control theory domain.
In \cite{GaLoRam}, Gasca, Lopez-Carmona, and Ramirez
proved a generalization of Sylvester's classical determinantal
identity; an alternative proof based on elimination techniques was
suggested in \cite{MuhlGas}, while a short and simpler one than the
earlier attempts was given in \cite{BeckMuhl}. More recently,
Mulders \cite{Mulders}  formulated another generalization related to
fraction free Gaussian elimination algorithms. Moreover, in order to
simplify the reading when dealing with minors identities, Leclerc
\cite{Leclerc} proposed a unified notation to encode many algebraic
relations between minors of a matrix, including Sylvester's
determinantal formula and the M\"uhlbach-Gasca generalization \cite{MuhlGas,Muhl}.

These generalizations (and sometimes only a rewriting of Sylvester's formula as Yakovlev's)
are scattered in the literature. Moreover, they used different notations which makes difficult to understand them as a whole and to see the connections between them. Therefore, one of the important aims of this paper is to present
these identities in a unified notation. Also, for helping the reader to better
understand these generalizations, they are illustrated by several examples.

Then, for showing the potentiality of our unified approach, a new generalization of Sylvester's identity will be given. Sylvester's determinantal identity can be recovered as a particular case of it. Such identities may have a potential
impact for the development of new recursive algorithms for interpolation and extrapolation.

In Section 2, the classical Sylvester's determinantal identity is
reminded, and some particular cases are discussed. Section 3 is
devoted to the generalizations, mentioned above, of the classical
identity. These identities are presented in a unified formulation
which allows to understand them better and to compare them.
In Section 4, the new generalization is proposed.

\section{The classical Sylvester's determinantal identity}
\label{SEC2}
 To introduce the classical Sylvester's determinantal identity we consider a square matrix $M=(a_{ij})$,
 of order $n$ over $\mathbb{K}$ (a commutative field). We denote its determinant by $\det M$
 and, for an integer
 $t$, $0\le t \le n-1$, and a given couple of integers  $i,j$ with $t < i,j \le n$,
 we define the
  determinant $a_{i,j}^{(t)}$ as

\begin{equation} \label{min} a_{i,j}^{(t)}= \left|
\begin{array}{ccc|c}
a_{11} & \cdots &a_{1t} & a_{1j} \\
a_{21} & \cdots &a_{2t} & a_{2j} \\
\vdots &  & \vdots  & \vdots \\
a_{t1} & \cdots &a_{tt} & a_{tj} \\
\hline
a_{i1} & \cdots &a_{it} & a_{ij} \\
\end{array}
\right|, \; {\rm for~} 0 < t \leq n-1 \;\;\; \mbox{and} \;\;\;
a_{i,j}^{(0)}=a_{ij},\end{equation} that is the determinant of order $(t + 1)$
obtained from the matrix $M$ by extending its leading principal
submatrix of order $t$ with the row $i$ and the column $j$ of $M$.
 Given these assumptions we can easily write down, in the following
Theorem, the Sylvester's identity.

\begin{theorem}[Sylvester's identity]
Let $M$ be a square matrix of order $n$ and $t$ an integer, $0 \le t
\le n-1$. Then, the following identity holds \begin{equation}  \label{IdSylv}
\det M \cdot \Big [ a_{t,t}^{(t-1)} \Big]^{n-t-1}= \left|
\begin{array}{ccc}
a_{t+1,t+1}^{(t)} & \cdots &  a_{t+1,n}^{(t)} \\
\vdots &  & \vdots \\
a_{n,t+1}^{(t)} & \cdots &  a_{n,n}^{(t)}
 \end{array}
\right|,
\end{equation}
with $a_{0,0}^{(-1)}=1$.
\end{theorem}
For different proofs of this Theorem, we refer to \cite{gant,Bair,AAM}.
We note also \cite{Konva}, where a combinatorial proof,  together with several non-commutative
extensions, has been presented.

Some particular cases of this Sylvester's identity are of interest:

\begin{description}
\item[\it Case 1:]
For $t=1$ and $a_{11} \not = 0$ we obtain the {\it Chi\'o pivotal
condensation} method \cite{chio} which states
$$
\det  M  =\frac{\det  B}{(a_{11})^{n-2}}
$$
where $B=(b_{ij})$, is the square matrix of order
$(n-1)$ with entries
$$b_{ij}=a_{i+1,j+1}^{(1)} = \left| \begin{array}{cc}
a_{11} & a_{1,j+1} \\ a_{i+1,1} & a_{i+1,j+1}
\end{array}
\right|,
$$
for $i,j=1,\ldots,n-1$.

\item[\it Case 2:]
For $t=n-2$,  we consider $M$ partitioned like

\begin{equation} \label{Msplit}
M=\left( \begin{array}{ccc}
a & {\mathbf{b}}^{T} & e \\
{\mathbf{c}} & D & {\mathbf{f}} \\
g & {\mathbf{h}}^{T} & l
\end{array}
\right) \end{equation} where $a$, $e$, $g$, $l$ are scalars, ${\mathbf{b}},
{\mathbf{c}}, {\mathbf{f}}, {\mathbf{h}}$ are $(n-2)$-length column
vectors, and $D$ is a  square matrix of order $(n-2)$. It is easy
 to apply the Sylvester's identity to such a matrix. In fact, by
 interchanging rows and columns in order to put $D$ in
the upper left corner, leaves the sign of determinant unaltered, so
$$
\det  M=\left| \begin{array}{ccc}
D & {\mathbf{c}} & {\mathbf{f}} \\
{\mathbf{b}}^{T} & a & e \\
{\mathbf{h}}^{T} & g & l
\end{array}
\right|,
$$
and applying (\ref{IdSylv}) we have

$$
\det M  \cdot \det D = \left| \begin{array}{cccc}
\left| \begin{array}{cc} D & {\mathbf{c}}
 \\ {\mathbf{b}}^{T} & a
\end{array}
\right| &
\left| \begin{array}{cc}
D & {\mathbf{f}} \\ {\mathbf{b}}^{T} & e
\end{array}
\right| \\ \\
 \left| \begin{array}{cc}
D & {\mathbf{c}} \\ {\mathbf{h}}^{T} & g
\end{array}
\right| &
\left| \begin{array}{cc}
D & {\mathbf{f}} \\ {\mathbf{h}}^{T} & l
\end{array}
\right|
\end{array}
\right|
=  \left| \begin{array}{cccc}
\left| \begin{array}{cc}
a & {\mathbf{b}}^{T} \\ {\mathbf{c}} & D
\end{array}
\right| &
\left| \begin{array}{cc}
{\mathbf{b}}^{T} & e \\ D & {\mathbf{f}}
\end{array}
\right| \\ \\
 \left| \begin{array}{cc}
{\mathbf{c}} & D \\ g & {\mathbf{h}}^{T}
\end{array}
\right| &
\left| \begin{array}{cc}
D & {\mathbf{f}} \\ {\mathbf{h}}^{T} & l
\end{array}
\right|
\end{array}
\right|.
$$
By setting
\begin{equation}
\label{blocchi}
A'=\left( \begin{array}{cc}
a & {\mathbf{b}}^{T} \\ {\mathbf{c}} & D
\end{array}
\right), \;
B'=\left( \begin{array}{cc}
{\mathbf{b}}^{T} & e \\ D & {\mathbf{f}}
\end{array}
\right), \;
C'=
\left( \begin{array}{cc}
{\mathbf{c}} & D \\ g & {\mathbf{h}}^{T}
\end{array}
\right), \;
D'=
\left( \begin{array}{cc}
D & {\mathbf{f}} \\ {\mathbf{h}}^{T} & l
\end{array}
\right),
\end{equation}
we can write the very simple rule
\begin{equation}
\label{SQUARE}
\det  M  \, \det  D =\det  A' \, \det  D' -\det  B' \,
\det  C',
\end{equation}
frequently used to obtain recursive algorithms in sequence
transformations (see, for instance, \cite{BR,BR1,Br-Rz}).
\end{description}

\section{Generalizations of  Sylvester's identity}

Several authors have deepened the main property of classical
Sylvester's identity, which is very useful in several domains.
 Some of these authors
have generalized the classical identity obtaining interesting new
formulae. Yakovlev  \cite{Yak1,Yak2} proved a generalization useful
in control theory. After,  Beckermann, Gasca, M\"uhlbach et al.
\cite{BeckMuhl,GaLoRam,MuhlGas} devoted some papers to generalizations
and applications of Sylvester's identity, while Mulders
\cite{Mulders} proposed a particular generalization of practical use
in linear programming and related topics. In this section, we
present a brief overview of these results.

In order to have a common working framework, we first introduce here
some general notations and definitions related to a $n \times m$ matrix
$M$ over $\mathbb{K}$.

 An {\it index list} $I$ of length $k \leq n$
is a $k$-tuple (with possible repetition) of integers
  from $\{1,2, \ldots, n\}$, that is
$I=(i_1,\ldots,i_k)$, with
$1 \leq i_l \leq n$ for $l=1,\ldots k$.
If we impose that all the elements are different, we have
an {\it index list without repetition}.
So, in this case,  $\text{card}(I) = k$. If, in addition,  $\alpha < \beta$
implies $i_\alpha< i_\beta $, then  $I$ is called an {\it ordered index list}.
The usual symbols $\in , \subset , \subseteq , \cup , \cap , \backslash,$ used for the sets
will be used for the ordered index lists and the result will be again an
ordered index list.
For any $n \in \mathbb{N}^+ $ (the set
of all positive integer), we define the ordered index list
$N_n=(1,2,\ldots,n)$.

In the sequel, if nothing is stated to the contrary, all index lists
are assumed to be ordered and without repetition.

Let now $I=(i_1,\ldots,i_{\alpha}) \subseteq N_n$ and
$J=(j_1,\ldots,j_{\beta}) \subseteq N_m$ two ordered index lists. We
denote the $ {\alpha} \times {\beta}$ submatrix, extracted from $M$,
with rows labeled by  $i \in I$ and columns labeled by  $j \in J$ as

\begin{equation}
\label{eqdefM}
M \left(  \!\!\! \begin{array}{c} I \\ J
 \end{array} \!\!\! \right) =
M \left( \!\!\!  \begin{array}{c} i_1, \ldots, i_{\alpha} \\ j_1,
\ldots, j_{\beta}
 \end{array}  \!\!\! \right)=
 \left( \begin{array}{ccc}
a_{i_1j_1} & \cdots &  a_{i_1j_{\beta}} \\
\vdots &  & \vdots \\
a_{i_{\alpha}j_1} & \cdots &  a_{i_{\alpha}j_{\beta}}
 \end{array}
\right).
\end{equation}
Let us mention that, in some papers, the indices for rows and columns are inverted. However, the notations we are using, seem to be the most common ones; see, for example, \cite{Aitk}.

If ${\alpha}={\beta}$, $n=m$, we denote the corresponding determinant as
$$
M \left[  \!\!\! \begin{array}{c} I \\ J
 \end{array} \!\!\! \right] =M \left[ \!\!\!  \begin{array}{c} i_1, \ldots, i_{\alpha} \\ j_1,
\ldots, j_{\alpha}
 \end{array}  \!\!\! \right]=
 \left| \begin{array}{ccc}
a_{i_1j_1} & \cdots &  a_{i_1j_{\alpha}} \\
\vdots &  & \vdots \\
a_{i_{\alpha}j_1} & \cdots &  a_{i_{\alpha}j_{\alpha}}
 \end{array}
\right|.
$$

With these notations, $a_{i,j}^{(t)} = M \left[ \!\!\!
\begin{array}{c} 1, \ldots, t, i \\ 1, \ldots, t, j
 \end{array}  \!\!\! \right]$, for $t< i,j \leq n$ and $0<t\leq n-1$.


\subsection{Generalization of Yakovlev}

In \cite{Yak1,Yak2}, Yakovlev presented some determinant identities
used to simplify several procedures in control theory. In
particular, in \cite{Yak2} he was the first to give  a simple generalization of
Sylvester's identity, which can be viwed as a renaming of the row and column indices.

Assume that $M=(a_{ij})$ is a square matrix of order $n$. Let $0<t\leq n-1$,
and consider the two ordered index lists $I=(i_1,\ldots,i_{t}) \subset N_n$, $J=(j_1,\ldots,j_{t})
\subset N_n$, and the complementary ordered index lists $I'=(i'_{t+1},\ldots,i'_{n}) \subset N_n$,
$J'=(j'_{t+1},\ldots,j'_{n}) \subset N_n$ (that is
 $I \cup I'=J \cup J'=N_n $).

Starting from the well-known Laplace expansion formula for the determinant of a squared matrix,
as proposed in \cite{gant},
the author proved a determinantal identity and use it in the solution of certain applied problems as
regulator synthesis for system of hight order.
This identity can be transformed in the following formula

\begin{equation} \label{yak}
\det  M \cdot \left( M \left[ \!\!  \begin{array}{l} i_1, \ldots,
i_t
\\ j_1, \ldots, j_t \end{array}\!\!  \right] \right) ^{n-t-1}=
\sum_{P(\alpha_{t+1}, \ldots,\alpha_{n})} \!\!(-1)^\mu \!\!
\prod_{\beta=t+1}^{n} M
\left[ \!\!\!  \begin{array}{l} i_1, \ldots, i_t,i'_{\alpha_\beta} \\
j_1, \ldots, j_t,j'_\beta \end{array}\!\!\! \right],
\end{equation}
where $P(\alpha_{t+1}, \ldots,\alpha_{n})$ represent the set of all
permutations of $( t+1,\ldots,n)$
and $\mu$ is the number of inversions needed to pass from $(t+1,
\ldots, n)$ to a certain permutation $(\alpha_{t+1}, \ldots, \alpha_n)$.

Let us now take $I=J=(1,\ldots,t)$. Thus,  $I'=J'=(t+1,\ldots,n)$. Formula
(\ref{yak}) becomes
$$
\det  M \cdot \left( M \left[ \!\!\!  \begin{array}{l} 1, \ldots, t
\\ 1, \ldots, t \end{array}\!\!\!  \right]  \right)  ^{n-t-1}=
\sum_{P(\alpha_{t+1}, \ldots,\alpha_n)}  \!\!(-1)^\mu \!\!
\prod_{\beta=t+1}^{n} M \left[ \!\!\!  \begin{array}{l} 1, \ldots,
t,\alpha_\beta \\ 1, \ldots, t,\beta \end{array} \!\!\! \right].
$$

This expression is one of the different ways to formulate the Sylvester's
determinantal identity.
In fact
$$M \left[ \!\!\!  \begin{array}{l} 1, \ldots,
t,\alpha_\beta \\ 1, \ldots, t,\beta \end{array} \!\!\! \right]=
a_{\alpha_{\beta},\beta}^{(t)}, $$
and, by
Leibniz formula,
\begin{eqnarray*}
\sum_{P(\alpha_{t+1}, \ldots,\alpha_n)}  \!\!(-1)^\mu \!\!
\prod_{\beta=t+1}^{n} a_{\alpha_{\beta},\beta}^{(t)}
\end{eqnarray*}
is exactly
the determinant of order $n-t$ on the right hand side of  (\ref{IdSylv}).

In the following examples we consider ${P(\alpha_{t+1}, \ldots,\alpha_n)}$ as an ordered set of $(n-t)!$ elements,
and we denote by $\mbox{\boldmath$\mu$}=(\mu_{1},\ldots, \mu_{(n-t)!})$ the vector whose elements are the
corresponding number of inversions.

\noindent{\bf Example 1:}

\noindent Let $n=6$ the order of $M$, and set $t=4$. We consider
the ordered index lists $I=(i_1,\ldots,i_4)=(1,3,5,6)$, $J=(j_1,\ldots,j_4)=(1,2,4,6)$;
then $I'=(i'_{5},i'_{6})=(2,4)$, $J'=(j'_{5},j'_{6})=(3,5)$. We have
$$
\det  M \cdot  M \! \left[ \!\!\!  \begin{array}{l} I
\\ J \end{array} \!\!\!\right]  =
\sum_{P(\alpha_{5},\alpha_6)}  (-1)^{\mbox{\footnotesize \boldmath$\mu$}} \prod_{\beta=5}^{6} \, M \!\left[ \!\!\!
\begin{array}{l} I,i'_{\alpha_\beta} \\ J,j'_\beta
\end{array}\!\!\! \right].
$$
Since $P(\alpha_{5}, \alpha_6)=\{(5,6),(6,5)\}$, and $\mbox{\boldmath$\mu$}=(0,1)$,  we obtain
$$
\det  M \cdot \! M \!\left[ \!\!\!  \begin{array}{l} I \\ J \end{array}\!\!\! \right]   =
M \left[ \!\!\!
\begin{array}{c} I,2 \\ J,3 \end{array}\!\!\!
\right] \cdot M \left[ \!\!\!  \begin{array}{c} I,4 \\ J,5
\end{array} \!\!\!\right]- M \left[ \!\!\!  \begin{array}{l} I,4 \\ J,3 \end{array} \!\!\!\right] \cdot M \left[
\!\!\!
\begin{array}{c} I,2 \\ J,5 \end{array}
\!\!\!\right].
$$

With the same values, but with $I=J=(1,\ldots,4)$, $I'=J'=(5,6)$,
$P(\alpha_{5}, \alpha_6)=\{(5,6),(6,5)\}$, and $\mbox{\boldmath$\mu$}=(0,1)$ we obtain
$$
\det  M \cdot  M \left[ \!\!\!  \begin{array}{l} I \\ J \end{array}\!\!\! \right]   = M \left[ \!\!\!
\begin{array}{c} I,5 \\ J,5 \end{array}\!\!\!
\right] \cdot M \left[ \!\!\!  \begin{array}{c} I,6 \\ J,6
\end{array} \!\!\!\right]- M \left[ \!\!\!  \begin{array}{l} I,6 \\ J,5 \end{array} \!\!\!\right] \cdot M \left[
\!\!\!
\begin{array}{c} I,5 \\ J,6 \end{array}
\!\!\!\right],
$$
that is
$$
\det  M \cdot  M \left[ \!\!\!  \begin{array}{c} 1, \ldots, 4 \\ 1,
\ldots, 4 \end{array}\!\!\! \right]= \left|
\begin{array}{cc}
a_{5,5}^{(4)} &   a_{5,6}^{(4)} \\
a_{6,5}^{(4)} &   a_{6,6}^{(4)}
 \end{array}
\right|,
$$
which is the classical Sylvester's identity.

\noindent{\bf Example 2:}

\noindent Let now  $n=8$ and $t=5$, and consider
the ordered index lists $I=(i_1,\ldots,i_5)=(2,3,5,6,8)$, $J=(j_1,\ldots,j_5)=(2,4,5,6,7)$, and the
complementary  $I'=(i'_{6},i'_{7},i'_{8})=(1,4,7)$ $J'=(j'_{6},j'_{7},j'_{8})=(1,3,8)$. We have
$$
\det  M \cdot \left( M \left[ \!\!\!  \begin{array}{l} I
\\ J \end{array} \!\!\!\right]  \right)  ^{2}=
\sum_{P(\alpha_{6},\alpha_7,\alpha_8)}  (-1)^{\mbox{\footnotesize \boldmath$\mu$}} \prod_{\beta=6}^{8} M \left[ \!\!\!
\begin{array}{l} I,i'_{\alpha_\beta} \\ J,j'_\beta
\end{array}\!\!\! \right].
$$
Since $P(\alpha_{6}, \alpha_7,\alpha_8)=
\{(8,7,6),(8,6,7),(7,8,6),(7,6,8),(6,7,8),(6,8,7)\}$, and $\mbox{\boldmath$\mu$}=(3,2,2,1,0,1)$,
we have
\begin{eqnarray*}
\det  M \cdot \left( M \left[ \!\!\!  \begin{array}{l} I\\ J \end{array}\!\!\! \right] \right)^2
\!\!\!& \!\!=\!\!& \!\!\!
-M \left[ \!\!\!
\begin{array}{c} I,7 \\ J,1 \end{array}\!\!\! \right] \cdot
M \left[ \!\!\!  \begin{array}{c} I,4 \\ J,3 \end{array} \!\!\!\right]
\cdot
M \left[ \!\!\!  \begin{array}{c} I,1 \\ J,8 \end{array} \!\!\!\right]+
M \left[ \!\!\!
\begin{array}{c}I,7 \\ J,1 \end{array}\!\!\! \right] \cdot
M \left[ \!\!\!  \begin{array}{c} I,1 \\ J,3 \end{array} \!\!\!\right]
\cdot
M \left[ \!\!\!  \begin{array}{c} I,4 \\ J,8 \end{array} \!\!\!\right]+   \\
&& \!\!\!\phantom{-}M \left[ \!\!\!
\begin{array}{c} I,4 \\ J,1 \end{array}\!\!\! \right] \cdot
M \left[ \!\!\!  \begin{array}{c} I,7 \\ J,3 \end{array} \!\!\!\right]
\cdot
M \left[ \!\!\!  \begin{array}{c} I,1 \\ J,8 \end{array} \!\!\!\right]-
M \left[ \!\!\!
\begin{array}{c} I,4 \\ J,1 \end{array}\!\!\! \right] \cdot
M \left[ \!\!\!  \begin{array}{c} I,1 \\ J,3 \end{array} \!\!\!\right]
\cdot
M \left[ \!\!\!  \begin{array}{c} I,7\\ J,8 \end{array} \!\!\!\right]+  \\
&& \!\!\!\phantom{-}M \left[ \!\!\!
\begin{array}{c} I,1 \\ J,1 \end{array}\!\!\! \right] \cdot
M \left[ \!\!\!  \begin{array}{c} I,4 \\ J,3 \end{array} \!\!\!\right]
\cdot
M \left[ \!\!\!  \begin{array}{c} I,7 \\ J,8 \end{array} \!\!\!\right]-
M \left[ \!\!\!
\begin{array}{c} I,1 \\ J,1 \end{array}\!\!\! \right] \cdot
M \left[ \!\!\!  \begin{array}{c} I,7 \\ J,3 \end{array} \!\!\!\right]
\cdot
M \left[ \!\!\!  \begin{array}{c} I,4 \\ J,8 \end{array} \!\!\!\right].
\end{eqnarray*}


\subsection{ Generalization of Gasca, Lopez-Carmona, and Ramirez}

In \cite{GaLoRam}, Gasca, Lopez-Carmona, and Ramirez proved a first generalization of
Sylvester's  determinantal identity; this generalized identity has
been applied to the derivation of a recurrence interpolation formula
for the solution of a general interpolation problem.
 Assume
 that $M$ is a square matrix of order $n$, and that
 $n=t+q$, for some positive integers $t,q$.
Given a set of $q$ ordered index lists $
J_k=(j_1^k,\ldots,j_{t+1}^k) \subset N_n$, $ k=1,\ldots,q,$ with
$$
\left\{
\begin{array}{ll}
\text{card} ( J_k)=t+1, & \quad k=1,\ldots,q,\\
\text{card} (J_k \cap J_{k+1})=t,  & \quad k=1,\ldots,q-1,
\end{array}
\right.
$$
we set
$$
\begin{array}{ll}
 \displaystyle J^{(q)}=\bigcup_{k=1}^{q} J_k \\
 S_k=J_k \cap J_{k+1}=(s_1^k,\ldots,s_{t}^k), \quad k=1,\ldots,q-1.
\end{array}
$$

Let $B$ the matrix with elements
$$
b_{ik}=M \left[ \!\!\!  \begin{array}{c} 1,2,
 \ldots, t, t+i\\
\!\!\! J_k
 \end{array}  \!\!\! \right] \qquad 1 \le i,k\le q
$$
and let $j_{h_k}^k$ be the element of $J_k$ such that
$$
\left\{
\begin{array}{ll}
j_{h_k}^k \in J_k-J_{k-1}=J_k-S_{k-1}, & \quad k=2,\ldots,q \\
\\
j_{h_1}^1 \in J_1-J_{2}=J_1-S_{1},  & \quad k=1.
\end{array}
\right.
$$
With these notations the authors give the following generalization of Sylvester's identity

\begin{equation} \label{Lop} \det  B= c \cdot \det  M \cdot
\prod_{k=1}^{q-1} M \left[ \!\!\!  \begin{array}{c} 1, \ldots, t\\
S_k
 \end{array}  \!\!\! \right]
\end{equation}

where $c$ is a sign factor which does not depend on the element
$a_{ij}$ of $M$, but only on the set of $J_k$ according to

$$
c= \left\{
\begin{array}{ll}
0 & \text{if card} (J^{(q)})< t+q \\
(-1)^{\mu}, \, \text{with} \, \mu=q(q-1)/2 + \displaystyle
\sum_{k=1}^{q}(j_{h_k}^k-h_k)  \,
   & \text{if card} (J^{(q)}) =t+q
\end{array}
\right.
$$
with $\displaystyle J^{(q)}=\bigcup_{k=1}^{q} J_k$.

  When $J_k=(1,\ldots,t,t+k)$, for
$k=1,2,\ldots,q$, we have, for all $k$,  $S_k=(1,\ldots,t)$, and the identity (\ref{Lop}) gives
$$
\det  B=  \det  M \cdot \left(   M \left[ \!\!\!
\begin{array}{c} 1,\ldots,t \\ 1,\ldots,t
 \end{array}  \!\!\! \right] \right)^{q-1}
 $$
which is exactly (\ref{IdSylv}), when $n=t+q$.

  An alternative proof of this identity was suggested by Gasca and
M\"uhlbach in \cite{MuhlGas}; this proof is based on the Laplacian
expansion and on an elimination strategy. Roughly speaking, when a
determinant is simplified by performing elementary operations
leaving its value unchanged, then {\it elimination strategies} are
applied, for example Gauss or Neville elimination.

\subsubsection{Examples}
\noindent{\bf Example 1:}

\noindent Let $M$ be a nonsingular square matrix of order $n=5$, choose $t=2, q=3$, and consider the ordered index lists
$J_1=(1, 3 ,4), J_2=(1, 4 ,5), J_3=(2 ,4 ,5)$ and $J^{(3)}=(1,2,3,4,5)=N_5$. We have $S_1=(1,4)$
and $S_2=(4 ,5)$. Let
$$
B = \left(
\begin{array}{ccc}
M\left[ \!\!\!\begin{array}{l} 1,2,3\\ 1,3,4 \end{array}
\!\!\!\right] &
M\left[ \!\!\!\begin{array}{l} 1,2,3\\ 1,4,5
\end{array} \!\!\!\right] &
M\left[ \!\!\!\begin{array}{l} 1,2,3\\ 2,4,5 \end{array} \!\!\!\right] \\
M\left[ \!\!\!\begin{array}{l} 1,2,4\\ 1,3,4 \end{array} \!\!\!\right] &
M\left[ \!\!\!\begin{array}{l} 1,2,4\\ 1,4,5 \end{array} \!\!\!\right] &
M\left[ \!\!\!\begin{array}{l} 1,2,4\\ 2,4,5 \end{array} \!\!\!\right] \\
M\left[ \!\!\!\begin{array}{l} 1,2,5\\ 1,3,4 \end{array} \!\!\!\right] &
M\left[ \!\!\!\begin{array}{l} 1,2,5\\ 1,4,5 \end{array} \!\!\!\right] &
M\left[ \!\!\!\begin{array}{l} 1,2,5\\ 2,4,5 \end{array} \!\!\!\right]
\end{array} \right).
$$
Since ${\rm card} (J^{(q)})=n$, for computing the sign factor we first have to determine $j^k_{h_k}$ for $k=1,\ldots,3$. For $k=1$, $j^1_{h_1}=3$,
$h_1=2$, for $k=2$, $j^2_{h_2}=5$, $h_2=3$, for $k=3$, $j^3_{h_3}=2$, $h_3=1$. So $\mu=3+4=7$ and $c=-1$.

\noindent Finally we have
$$
\det B = - \det M \cdot M \left[ \!\!\!\begin{array}{l} 1,2\\ 1,4 \end{array} \!\!\!\right]
\cdot M \left[ \!\!\!\begin{array}{l} 1,2\\ 4,5 \end{array} \!\!\!\right].
$$

 \noindent{\bf Example 2:}
\noindent Let again $n=5, t=2, q=3$, but consider the ordered index lists
$J_1=(1 ,2, 3), J_2=(2, 3, 4), J_3=(1, 2 ,4)$ with ${\rm card}(J^{(3)}=(1,2,3,4))<5$. So $c=0$, and the
determinant of the matrix
$$
B = \left(
\begin{array}{ccc}
M\left[ \!\!\!\begin{array}{l} 1,2,3\\ 1,2,3 \end{array}
\!\!\!\right] &
M\left[ \!\!\!\begin{array}{l} 1,2,3\\ 2,3,4
\end{array} \!\!\!\right] &
M\left[ \!\!\!\begin{array}{l} 1,2,3\\ 1,2,4 \end{array} \!\!\!\right] \\

M\left[ \!\!\!\begin{array}{l} 1,2,4\\ 1,2,3 \end{array} \!\!\!\right] &
M\left[ \!\!\!\begin{array}{l} 1,2,4\\ 2,3,4 \end{array} \!\!\!\right] &
M\left[ \!\!\!\begin{array}{l} 1,2,4\\ 1,2,4 \end{array} \!\!\!\right] \\

M\left[ \!\!\!\begin{array}{l} 1,2,5\\ 1,2,3 \end{array} \!\!\!\right] &
M\left[ \!\!\!\begin{array}{l} 1,2,5\\ 2,3,4 \end{array} \!\!\!\right] &
M\left[ \!\!\!\begin{array}{l} 1,2,5\\ 1,2,4 \end{array} \!\!\!\right]
\end{array} \right)
$$
 is equal to zero.


\subsection{Generalization of Beckermann, Gasca and M\"uhlbach}

Sylvester's classical identity can be interpreted as an extension of
Leibniz's definition of a matrix determinant,  known as the {\it
Muir's law of extensible minors} \cite{Muir}. Muir
first stated the {\it Law of Complementaries}, which was been
already known to Cayley in 1878 and reads as

{\it To every general theorem which takes the form of an identical relation between
a number of the minors of a determinant or between the determinant itself and
a number of its minors, there corresponds another theorem derivable from the
former by merely substituting for every minor its cofactor in the determinant,
and then multiplying any term by such a power of the determinant  that will make
the terms of the same degree.}

Brualdi and Schneider \cite{BrSch} gave a formal treatment of
determinantal identities of the minors of a matrix and then provided
a careful exposition of the {\it Law of Extensible Minors} and of
the {\it Law of Complementaries}, two methods for obtaining from a
given determinantal identity another determinantal identity, as
mentioned earlier. In  \cite{Muhl}, M\"uhlbach using the relation
between Muir's Law of Extensible Minors, Sylvester's identity and
the Schur complement, presented a new principle for extending
determinantal identities which generalizes Muir's Law. M\"uhlbach's
proof makes use of general elimination strategies and of generalized
Schur complements.  As applications of this technique, he derived a
generalization of Sylvester's identity, which is the same as the one proposed
by Gasca and M\"uhlbach in \cite{MuhlGas}. Later, Beckermann and
M\"uhlbach \cite{BeckMuhl} gave an approach, shorter and
conceptually simpler than the earlier attempts \cite{MuhlGas,Muhl},
for a general determinantal identity of Sylvester's type. In the
sequel we report these last results.

  We consider a  $n
\times m$ matrix $M$ with rows and columns numbered in the usual way.
Let $I=(i_1,\ldots,i_{\alpha}) \subset N_n$ and
$J=(j_1,\ldots,j_{\beta}) \subset N_m$ two ordered index lists, and
consider also the row vectors $
{\mathbf{z}}_1,\ldots,{\mathbf{z}}_q$, where
${\mathbf{z}}_k=(z_{k,j})_{j \in N_m}$.

 We border the matrix $M \left( \!\!\!  \begin{array}{c} I\\ J
 \end{array}  \!\!\! \right)$,
extracted from $M$, by the $q$ row vectors
${\mathbf{z}}^{\prime}_k=(z_{k,{j_{\lambda}}})_{\lambda=1,\ldots,{\beta}}$,
$k=1,\ldots,q$, extracted from the ${\mathbf{z}}_k$, that is we
consider
$$
M
\left( \!\!\!  \begin{array}{c}
i_1, \ldots ,i_{\alpha} \vert \, 1, \ldots , q \\
j_1,  \ldots, j_{\beta}
 \end{array}  \!\!\! \right)=
 \left( \begin{array}{ccc}
 a_{i_1,j_1} & \cdots & a_{i_1,j_{\beta}} \\
 \vdots & & \vdots \\
 a_{i_{\alpha},j_1} & \cdots & a_{i_{\alpha},j_{\beta}} \\
 \\
 \hline
z_{1,{j_1}} &\cdots & z_{1,{j_{\beta}}} \\
\vdots & & \vdots \\
z_{q,{j_1}} &\cdots & z_{q,{j_{\beta}}} \\
\end{array} \right).
 $$

Let $q \in  \mathbb{N}^+$ be fixed, and consider the index lists
$$
I_1,\ldots, I_q \subset N_n \quad \text{and} \quad J_1,\ldots, J_q
\subset N_m
$$
such that
$$
\text{\rm card}(J_k)=\text{\rm card}(I_k)+1 \quad \text{for} \quad
k=1,\ldots,q.
$$
We set, for $k=1,\ldots,q$,
$$
I^{(k)}:= \bigcap_{i=1}^{k} I_i \quad \text{and} \quad J^{(k)}:=
\bigcup_{j=1}^{k} J_j
$$

The authors give the following theorem.

\begin{theorem}
\label{TH1}
Let
$$
I_0 \subseteq I^{(q)} \mbox{~~such that~~} \text{\rm
card}(I_0)=\text{\rm card}(J^{(q)})-q.
$$
Then, for any matrix $N_n \times N_m$ matrix $M$, there exists an
element $c \in \mathbb{K}$ depending only on $M \left( \!\!\!
\begin{array}{c} I_0 \\ J^{(q)}
 \end{array}  \!\!\! \right)$ such that,
from all row vectors
 ${\mathbf{z}}_1,\ldots,{\mathbf{z}}_q \in \mathbb{K}^{m}$, and so
 for the vectors
${\mathbf{z}}^{\prime}_k=(z_{k,{j_{\lambda}}})_{{j_{\lambda}}\in J^{(q)}}$, $k=1, \ldots, q$,
extracted from them,

\begin{equation} \label{BecMul} \det  B= c \cdot   M \left[ \!\!\!
\begin{array}{c}
 I_0 \vert \, 1, \ldots, q \\ J^{(q)}
 \end{array}  \!\!\! \right],
\end{equation}
where $B$ is the matrix with elements
$$
b_{ij}= M \left[ \!\!  \begin{array}{c}
 I_i \vert \, j \\ J_i
 \end{array}  \!\! \right], \qquad 1 \le i,j\le q,
$$
that is the matrix $M \left[ \!\!  \begin{array}{c}
 I_i \\ J_i
 \end{array}  \!\! \right]$ bordered with the row vector
${\mathbf{z}}^{\prime}_j=(z_{j,{j_{\lambda}}})_{\lambda\in J_i}$.
\end{theorem}
Two corollaries in \cite{BeckMuhl} are also of interest.
\begin{corollary}
\label{COR1}
If ${\rm card}(I^{(k)}) > {\rm card}(J^{(k)})-k$, for some $k \in \{2, \ldots , q\}$,
then (\ref{BecMul}) holds with $c=0$.
\end{corollary}
\begin{corollary}
\label{COR2}
If ${\rm card}(J^{(q)}) <q$, then, for all choices of vectors
 ${\mathbf{z}}_1,\ldots,{\mathbf{z}}_q$, we have $\det B=0$.
\end{corollary}

 In some particular cases, the constant $c$ of Theorem \ref{TH1} can be computed explicitly. This is
the case of Sylvester's classical determinantal identity. In fact, we consider a nonsingular
matrix $M$ of order
$n$,  and we choose two integer $t,q \in \mathbb{N}^+$, such that
$t+q=n$. Let
$$
\begin{array}{lr}
J_k=(1,\ldots,t,t+k), & \text{for} \quad k =1,\ldots,q \\
\\
I_k=(1,\ldots,t), & \text{for} \quad k =1,\ldots,q \\
\\
I_0=I^{(q)}=(1,\ldots,t), \quad J^{(q)}=N_n=(1,\ldots,n),& \\
\\
{\mathbf{z}}^\prime_k=M \left( \!\!\!  \begin{array}{c} t+k \\ J_k
 \end{array}  \!\!\! \right), &  \text{for} \quad k =1,\ldots,q.
\end{array}
$$
We consider the matrix  $B$  with elements
$$
b_{ij}=
M \left[ \!\!\! \begin{array}{c}  1,\ldots,t \vert \, i\\ 1,\ldots,t, t+j
 \end{array}  \!\!\! \right],
\qquad 1 \le i,j\le q.
 $$
  From (\ref{BecMul}) we have
\begin{equation} \label{BSylv} \det  B= c \cdot   M \left[ \!\!\!
\begin{array}{c} 1,\ldots,t+q \\ 1,\ldots,t+q
 \end{array}  \!\!\! \right] = c \cdot \det M
\end{equation}
where $c$ depends only on $ M \left( \!\!\!  \begin{array}{c}
 I_0 \\ N_m \end{array}  \!\!\! \right)$ and can be evaluated by choosing the particular row vectors
${\mathbf{z}}_k=(\delta_{\lambda,t+k})_{\lambda \in N_n}$ for
$k=1,\ldots,q$, where $\delta$ is the usual Kronecker symbol. With this
choice we have
$$
\det  B =\left(   M \left[ \!\!\!  \begin{array}{c} 1,\ldots,t \\
1,\ldots,t
 \end{array}  \!\!\! \right] \right)^q,
 $$
while the right hand side becomes
$$
c \cdot   M \left[ \!\!\!  \begin{array}{c} 1,\ldots,t \\
1,\ldots,t
 \end{array}  \!\!\! \right].
$$
So, assuming that $\displaystyle M \left[\!\!\!  \begin{array}{c} 1,\ldots,t \\
1,\ldots,t
 \end{array}  \!\!\! \right] \not = 0$,
$$
c=\left(  M \left[\!\!\!  \begin{array}{c} 1,\ldots,t \\
1,\ldots,t
 \end{array}  \!\!\! \right] \right)^{q-1},
$$
and from (\ref{BSylv}) it follows
$$
\det  B=  \det  M \cdot \left(   M \left[ \!\!\!
\begin{array}{c} 1,\ldots,t \\ 1,\ldots,t
 \end{array}  \!\!\! \right] \right)^{q-1},
 $$
which is the classical Sylvester's identity (\ref{IdSylv}).

 In a similar way \cite{BeckMuhl},  is also possible to obtain  the identities of Schweins and
Monge \cite{Aitk}.

The relation (\ref{BecMul}) can be very useful for obtaining
other determinantal identities, since $c$ depends on $M$ but not on
the rows corresponding to indexes from $N_n \setminus I_0$.
As an application of this general determinantal
identity, in \cite{BeckMuhl} old and new
 recurrence relations for the E-transforms \cite{BR} are derived. In the same paper, a more general Theorem,
 based on a {\it chain} of index lists (that is lists $I_k\subset N_n$ and $J_k \subset N_m$ such that, for all
 $k=1, \ldots, q$, ${\rm card}(J^{(k)})={\rm card}(I^{(k)})+k$) is also provided.

\subsubsection{Examples}
\noindent{\bf Example 1:}

\noindent Let $M$ a $n \times m$ matrix, with $n=6$ and $ m=7$. Set $q=3$, and consider the ordered index lists
$I_1=(2, 3, 4),
I_2=(2 , 4),
I_3=(1, 2 ,4)$ and
$J_1=(2, 3, 4 ,5),
J_2=(2 ,3 ,4),
J_3=(2, 3 ,4 ,7)$.  We have $I^{(3)}=(2,4)$
and $J^{(3)}=(2   ,  3  ,   4   ,  5  ,   7)$. Choose any ${\mathbf{z}}_1, {\mathbf{z}}_2,
{\mathbf{z}}_3$, and $I_0=I^{(3)}$. Let
$$
B = \left(
\begin{array}{ccc}
M\left[ \!\!\!\begin{array}{l} 2,3,4 \,|\, 1\\ 2  ,   3   ,  4   ,  5 \end{array}
\!\!\!\right] &
M\left[ \!\!\!\begin{array}{l} 2,3,4\,|\,2\\ 2  ,   3   ,  4   ,  5
\end{array} \!\!\!\right] &
M\left[ \!\!\!\begin{array}{l} 2,3,4\,|\,3\\  2  ,   3   ,  4  ,5 \end{array} \!\!\!\right] \\
M\left[ \!\!\!\begin{array}{l} 2,4\,|\,1\\ 2  ,   3   ,  4   \end{array} \!\!\!\right] &
M\left[ \!\!\!\begin{array}{l} 2,4\,|\,2\\ 2  ,   3   ,  4   \end{array} \!\!\!\right] &
M\left[ \!\!\!\begin{array}{l} 2,4\,|\,3\\ 2  ,   3   ,  4   \end{array} \!\!\!\right] \\
M\left[ \!\!\!\begin{array}{l} 1,2,4\,|\,1\\  2   ,  3  ,   4   ,  7 \end{array} \!\!\!\right] &
M\left[ \!\!\!\begin{array}{l} 1,2,4\,|\,2\\ 2   ,  3  ,   4   ,  7 \end{array} \!\!\!\right] &
M\left[ \!\!\!\begin{array}{l} 1,2,4\,|\,3\\ 2   ,  3  ,   4   ,  7 \end{array} \!\!\!\right]
\end{array} \right),
$$
where $\det B \not = 0$.
 If $M\left[ \!\!\!\begin{array}{l} 2,4 \,|\, 1,2,3\\ 2   ,  3  ,   4   ,  5  ,   7 \end{array}
\!\!\!\right] \not = 0$, there exists $c$ (which cannot be computed a priori) such that (\ref{BecMul}) holds.

 \noindent{\bf Example 2:}
\noindent Let again $n=6$ and $m=7$, but consider $q=5$ and the ordered index lists
$I_1=(1,4),
I_2=(3 ,4),
I_3=(2, 4),
I_4=(4 , 5),
I_5=(4, 6)
$
and
$J_1=(1 ,2, 4),
J_2=(1, 4 ,7),
J_3=(1, 2 ,7),
J_4=(1, 4 ,7),
J_5=(1, 2 ,7)$.  We have $I^{(5)}=(4)$
and $J^{(5)}=(1   ,  2   ,  4   ,  7)$.  Let
$$
B = \left(
\begin{array}{ccccc}
M\left[ \!\!\!\begin{array}{l} 1,4 \,|\, 1\\ 1,2,4 \end{array} \!\!\!\right] &
M\left[ \!\!\!\begin{array}{l} 1,4 \,|\,2\\ 1,2,4 \end{array} \!\!\!\right] &
M\left[ \!\!\!\begin{array}{l} 1,4\,|\,3\\  1,2,4\end{array} \!\!\!\right] &
M\left[ \!\!\!\begin{array}{l} 1,4 \,|\,4\\ 1,2,4 \end{array} \!\!\!\right] &
M\left[ \!\!\!\begin{array}{l} 1,4 \,|\,5\\ 1,2,4 \end{array} \!\!\!\right] \\
M\left[ \!\!\!\begin{array}{l} 3,4\,|\,1\\ 1,4,7  \end{array} \!\!\!\right] &
M\left[ \!\!\!\begin{array}{l} 3,4\,|\,2\\ 1,4,7  \end{array} \!\!\!\right] &
M\left[ \!\!\!\begin{array}{l} 3,4\,|\,3\\ 1,4,7  \end{array} \!\!\!\right] &
M\left[ \!\!\!\begin{array}{l} 3,4\,|\,4\\ 1,4,7  \end{array} \!\!\!\right] &
M\left[ \!\!\!\begin{array}{l} 3,4\,|\,5\\ 1,4,7  \end{array} \!\!\!\right]  \\
M\left[ \!\!\!\begin{array}{l} 2,4\,|\,1\\  1,2,7 \end{array} \!\!\!\right] &
M\left[ \!\!\!\begin{array}{l} 2,4\,|\,2\\ 1,2,7  \end{array} \!\!\!\right] &
M\left[ \!\!\!\begin{array}{l} 2,4\,|\,3\\ 1,2,7  \end{array} \!\!\!\right] &
M\left[ \!\!\!\begin{array}{l} 2,4\,|\,4\\ 1,2,7  \end{array} \!\!\!\right] &
M\left[ \!\!\!\begin{array}{l} 2,4\,|\,5\\ 1,2,7  \end{array} \!\!\!\right]
\end{array} \right).
$$
Due to Corollary \ref{COR2}, since ${\rm card}(J^{(5)})<5$, then, for any ${\mathbf{z}}_1, \ldots,
{\mathbf{z}}_5$, $\det B$ is equal to zero.


\subsection{Generalization of { Mulders}}

Another generalization of Sylvester's identity is due to Mulders
\cite{Mulders}. The author follows the idea of Bareiss \cite{Bair,Bair2}, where the Sylvester's
identity is used to prove that certain
Gaussian elimination algorithms, transforming a matrix into an
upper-triangular form, are fraction-free.
Fraction-free algorithms also have applications in the computation of matrix rational approximants, matrix GCDs, and
generalized Richardson extrapolation processes \cite{BeckLaba}. They are used
for controlling, in exact arithmetic, the growth of intermediate results.
The generalization of Sylvester's identity, due to Mulders, was used to prove that also
certain {\it random} Gaussian elimination algorithms are fraction-free.
The fraction-free random Gaussian elimination algorithm, obtained  in that way, has been used
in the simplex method (to solve linear programming problems) and for finding a
solution of $A {\mathbf{x}} = {\mathbf{b}}, {\mathbf{x}}\geq 0$, with $A$ a $m \times n$ matrix.

 Given a $n \times m$ matrix $M$, we consider the
index list $I=(i_1,\ldots,i_t) \subseteq N_n$ (without repetition but not necessarily ordered)
and another index list $J=(j_1,\ldots,j_t)
\subseteq N_m$, (with possible repetition and not necessarily ordered). The two lists have the same length
$t$ ,
 $0 \leq t \leq \min(m,n)$. Notice that it is possible to have $t=0$, so an empty list.
With  these assumptions, we can define an equivalence class $[(i_1,j_1),\ldots,(i_t,j_t)]$ of all possible
permutations of the $t$ pairs and define the determinant
$$
 a^{[(i_1,j_1),\ldots,(i_t,j_t)]} =
 M \left[ \!\!\!  \begin{array}{c} i_1, \ldots, i_t \\ j_1,
\ldots, j_t
 \end{array}  \!\!\! \right],
$$
 with $a^{[\O]}=M\left[\begin{array}{c}
 \O \\ \O \end{array}\right]=1$ when $t=0$. Note that, when $J$ has repeated elements, the determinant is equal to
 zero.

We also define the operation $\leftarrow$ on the pairs defining the above determinant as
$$ [(i_1,j_1),\ldots,(i_t,j_t)] \leftarrow [(u_1,v_1),\ldots,(u_r,v_r)].$$
$I$ and $U=(u_1,\ldots,u_r)$ are index lists without repetition but not necessarily ordered,
$J$ and $V=(v_1,\ldots,v_r)$ are index lists not necessarily ordered but with possible repetitions.
This operation produces a new class of pairs by adding to the first class all the pairs $(u_{k}, v_{k})$
of the second class when $u_{k} \not \in I$, and by replacing the pair  $(i_{\ell}, j_{\ell})$ of the
first class by
the pair $(u_{k}, v_{k})$, when  $u_{k}=i_{\ell}$ for some $i_\ell \in I$. This operation is performed for $k=1,\ldots,r$,
where $r$ is an integer.
This operation is associative and commutative if
$I \cap U = \O$.
For example, the result of $ [(2,1),(3,3),(1,5),(5,3)] \leftarrow
[(1,1),(2,3),(4,4)],$ is  \([(2,3),(3,3),(1,1),(5,3),(4,4)]\),
since the pair $(4,4)$ has been added, and the pairs $(1,1), (2,3)$ have
replaced the pairs $(1,5), (2,1)$, respectively.

 We consider now two integers $i$ and $j$, $1 \leq i \leq n, 1 \leq j \leq
 m$ and the
 determinant
\begin{equation} \label{min2}
 a^{[(i_1,j_1),\ldots,(i_t,j_t)] \leftarrow [(i,j)]} .
\end{equation}
Due to the operation previously defined, the result is the following:

\begin{itemize}
\item When $i \not \in I$,  then we simply
extend the matrix with the $i$th row and $j$th column of $M$, and we
compute the determinant of the new square matrix of order $(t+1)$,
that is
$$a^{[(i_1,j_1),\ldots,(i_t,j_t)] \leftarrow [(i,j)]}=
M \left[ \!\!\!  \begin{array}{c}
i_1, \ldots, i_t,i \\ j_1, \ldots, j_t,j
 \end{array}  \!\!\! \right].$$

\item When it exists $1 \le k \le t$ such that  $i_k=i$,
 we replace, in $a^{[(i_1,j_1),\ldots,(i_t,j_t)]}$, the pair $(i_k,j_k)$ by the pair
 $(i,j)$, that is we take the following determinant of a new matrix of order $t$
 $$a^{[(i_1,j_1),\ldots,(i_t,j_t)] \leftarrow [(i,j)]}=
M \left[ \!\!\!  \begin{array}{c} i_1, \ldots,i_k\!\!=\!i\,, \ldots, i_t \\
j_1, \ldots,\phantom{i~}j\phantom{~~~}, \ldots, j_t
 \end{array}  \!\!\! \right] .$$

\end{itemize}

Let us consider now the  particular case where  $i_k=k, j_k=k$, for all
$k=1,\ldots,t$, that is $I=J=N_t$. When $t< i,j \leq \min (n,m)$, we recover the usual
definition $a_{i,j}^{(t)}$ given in (\ref{min}), but not in the other cases.
For simplicity, in the sequel, we set
$\widetilde{a}_{i,j}^{(t)}= a^{[(1,1),\ldots,(t,t)]\leftarrow [(i,j)]}$.
 We have
\begin{eqnarray*}
\left\{
\begin{array}{l@{~~~~}l@{~~~~~~~}l}
{\rm If~} t< i,j &  \widetilde{a}_{i,j}^{(t)} \mbox{~of order $t+1$} & \widetilde{a}_{i,j}^{(t)}={a}_{i,j}^{(t)}\\
{\rm If~} j \leq t < i &  \widetilde{a}_{i,j}^{(t)} \mbox{~of order $t+1$} & \widetilde{a}_{i,j}^{(t)}=0\\
{\rm If~} i,j \leq t, i \not = j& \widetilde{a}_{i,j}^{(t)} \mbox{~of order $t$} & \widetilde{a}_{i,j}^{(t)}=0\\
{\rm If~} i =j \leq t& \widetilde{a}_{i,j}^{(t)} \mbox{~of order $t$} &  \widetilde{a}_{i,j}^{(t)}=
a^{[(1,1),\ldots,(t,t)]}\\
{\rm If~} i \leq t<j & \widetilde{a}_{i,j}^{(t)} \mbox{~of order $t$} &
 \widetilde{a}_{i,j}^{(t)} \mbox{\rm ~is a new determinant.}
\end{array}
\right.
\end{eqnarray*}

By using this particular choice for $I$ and $J$,  Mulders
formulated the following generalized Sylvester's determinantal
identity that allows to consider the determinantal elements
$\widetilde{a}_{i,j}^{(t)}$ when $i$ and/or $j$ are $\le t$.
 In the sequel, we state a more generalized identity.

\begin{theorem}

For $0 \le t \le \min(n,m)$, $0 \le p,q \le t$ and $1 \le s \le
\min(n-p,m-q)$, the following identity holds
\begin{equation}
\label{Mulderseq}
a^{[(1,1),\ldots,(t,t)]\leftarrow [(p+1,q+1),\ldots,(p+s,q+s)]}
\Big[ \widetilde{a}_{t,t}^{(t-1)} \Big]^{s-1} =
\left| \begin{array}{ccc}
\widetilde{a}_{p+1,q+1}^{(t)} & \cdots &  \widetilde{a}_{p+1,q+s}^{(t)} \\
\vdots &  & \vdots \\
\widetilde{a}_{p+s,q+1}^{(t)} & \cdots &  \widetilde{a}_{p+s,q+s}^{(t)}
 \end{array}
\right|.
\end{equation}
\end {theorem}

For $m=n,p=q=t$ and $s=n-t$ we obtain the classical Sylvester's
identity (\ref{IdSylv}).

\subsubsection{Examples}
\noindent{\bf Example 1:}

\noindent Let $M$ a $n \times m$ matrix, with $n=7$ and $ m=8$. Set $t=5$, and $I=J=N_5$. Assume $p=3,q=4$,
and $s=3$, so in this case $p+s \geq t$. Then
we have
$\displaystyle a^{[(1,1),\ldots,(5,5)]\leftarrow [(4,5),(5,6),(6,7)]} =
M \left[ \!\!\!  \begin{array}{c} 1, 2,3, 4,5,  6 \\ 1,
2, 3, 5, 6, 7  \end{array}\!\!\! \right]$ and
$\displaystyle \widetilde{a}_{5,5}^{(4)} = M \left[ \!\!\!  \begin{array}{c} 1, \ldots,5 \\ 1,
\ldots,  5 \end{array}\!\!\! \right]$. The right hand side of (\ref{Mulderseq}) is the following
determinant of order 3
$$
\left| \begin{array}{ccc}
\widetilde{a}_{4,5}^{(5)} & \widetilde{a}_{4,6}^{(5)} &  \widetilde{a}_{4,7}^{(5)} \\
\widetilde{a}_{5,5}^{(5)} & \widetilde{a}_{5,6}^{(5)} &  \widetilde{a}_{5,7}^{(5)} \\
\widetilde{a}_{6,5}^{(5)} & \widetilde{a}_{6,6}^{(5)} &  \widetilde{a}_{6,7}^{(5)}
 \end{array}
\right|,
$$
that is, thanks to the operation defined,
$$
 \left|
\begin{array}{ccc}
0 &
M\left[ \!\!\!\begin{array}{l} 1,2,3,4,5\\ 1,2,3,6,5
\end{array} \!\!\!\right] &
M\left[ \!\!\!\begin{array}{l} 1,2,3,4,5\\ 1,2,3, 7,5 \end{array} \!\!\!\right] \\
M\left[ \!\!\!\begin{array}{l} 1,\ldots,5\\ 1,\ldots,5  \end{array} \!\!\!\right] &
M\left[ \!\!\!\begin{array}{l} 1,2,3,4,5\\ 1,2,3,4,6  \end{array} \!\!\!\right] &
M\left[ \!\!\!\begin{array}{l} 1,2,3,4,5\\ 1,2,3,4,7  \end{array} \!\!\!\right] \\
0 &
M\left[ \!\!\!\begin{array}{l}1,\ldots,6\\ 1,\ldots,6\end{array} \!\!\!\right] &
M\left[ \!\!\!\begin{array}{l} 1,2,3,4,5,6\\ 1,2,3,4,5,7 \end{array} \!\!\!\right]
\end{array} \right|,
$$
and (\ref{Mulderseq}) holds.

\noindent{\bf Example 2:}

\noindent Let again $M$ a $7 \times 8$ matrix. Set now $t=6$, and $I=J=N_6$. Assume $p=2,q=3$,
and $s=3$. In this case $p+s <t$, so
we have
$\displaystyle a^{[(1,1),\ldots,(6,6)]\leftarrow [(3,4),(4,5),(5,6)]} =
M \left[ \!\!\!  \begin{array}{c} 1, 2,3, 4,5,6\\ 1,
2, 4, 5, 6,6 \end{array}\!\!\! \right]$ and
$\displaystyle \widetilde{a}_{6,6}^{(5)} = M \left[ \!\!\!  \begin{array}{c} 1, \ldots,6 \\ 1,
\ldots,  6 \end{array}\!\!\! \right]$. It follows that the left hand side is equal to zero for the presence of two equal columns in the first determinant.
The right hand side  is now
$$
\left| \begin{array}{ccc}
\widetilde{a}_{3,4}^{(6)} & \widetilde{a}_{3,5}^{(6)} &  \widetilde{a}_{3,6}^{(6)} \\
\widetilde{a}_{4,4}^{(6)} & \widetilde{a}_{4,5}^{(6)} &  \widetilde{a}_{4,6}^{(6)} \\
\widetilde{a}_{5,4}^{(6)} & \widetilde{a}_{5,5}^{(6)} &  \widetilde{a}_{5,6}^{(6)}
 \end{array}
\right|=
 \left|
\begin{array}{ccc}
0 & 0 & 0\\
M\left[ \!\!\!\begin{array}{l} 1,\ldots,6\\ 1,\ldots,6  \end{array} \!\!\!\right] &
0 & 0 \\
0 &
M\left[ \!\!\!\begin{array}{l}1,\ldots,6\\ 1,\ldots,6\end{array} \!\!\!\right] &
0
\end{array} \right| =0,
$$
and (\ref{Mulderseq}) holds.

\section{A new generalization}

We present now a different generalization of Sylvester's identity.
Let $M$ be a matrix of dimension $n \times m$. Chosen two
indices $t$ and $s$,  $0 < t \le r$, with $r = \min(n,m)$, and $1 \le s \le \min(n-t,m-t)$, in such
a way that it exists $q \in \mathbb{N}$  which satisfies $r=t+q
\, s$,
we  can build a sequence of submatrices $M_{k}$ of $M$,

\begin{equation} \label{Mk}
M_{k}= M \left( \!\!\! \begin{array}{c}
1 ,\ldots,t, t+1, \ldots ,t+ks \\1 ,\ldots,t, t+1, \ldots ,t+ks \end{array} \!\!\! \right)\!, \, \,
\mbox{for~}
0\leq k \leq q
\end{equation}
where, by definition, $\det M_{0} = a_{t,t}^{(t-1)}$ and $\det M_{q} = a_{r,r}^{(r-1)}$.
We avoid to indicate the dependence from $t$ and $s$ since they are fixed.
So, for the same reason, we set $m_{k}=t+ks$, for $k=0,\ldots, q-1$, and
we consider the matrix $B_k$ of order $s$, with elements
\begin{equation} \label{BIJ}
b_{ij}^{(k)}= M \left[ \!\!  \begin{array}{c}
1, \ldots,  m_k, m_k+ i  \\ 1, \ldots,  m_k, m_k+ j
 \end{array}  \!\! \right] = a_{m_k+i,m_k+j}^{(m_k)} \,, \qquad 1 \le i,j\le s.
\end{equation}

We remind that, by definition,  when $\alpha > \beta$, the product $\displaystyle \prod_{i=\alpha}^\beta \odot $
is  empty, and its value is set to $1$, and that, when $0^0$ is considered
 as an empty product of zeros, its value is also set to $1$. The following
theorem holds.

\begin{theorem}[Generalized Sylvester's identity] \label{th}

Let $M$ be a matrix of dimension $n \times m$, and $M_0$ its square leading principal
submatrix of fixed order $t$, $0 < t \le \min(n,m)$. Then, chosen $s$, $q \in
\mathbb{N} $,  $1 \le s \le \min(n-t,m-t)$, such that
$ \min(n,m)=t+q \, s$, the following identities hold for $0 \leq  k \leq q$

\begin{equation} \label{GenSylvE} \frac{\det\, M_{k}} {\left[  \det\,
M_{0} \right] ^{(s-1)^k} }= \frac{\displaystyle
\prod_{\substack{ i=1 \\ i  \text{\rm ~odd}}}^{k-1}
 \,\left[  \det B_{i} \right]^{(s-1)^{k-1-i}}}
{\displaystyle \prod_{\substack{ j=0 \\ j   \text{\rm ~even}}}^{k-2}
 \,  \left[  \det B_{j} \right]^{(s-1)^{k-1-j}}}\,\,, \qquad \text{for} \,\, k  \text{\rm ~even},
\end{equation}

\begin{equation} \label{GenSylvO} \det\, M_{k} \, \cdot  \left[ \det\,
M_{0} \right] ^{(s-1)^k}= \frac{\displaystyle \prod_{\substack{
i=0 \\ i  \text{\rm ~even}}}^{k-1}
 \,\left[  \det B_{i} \right]^{(s-1)^{k-1-i}}}
{\displaystyle \prod_{\substack{ j=1 \\ j  \text{\rm ~odd}}}^{k-2}
 \,\left[  \det B_{j} \right]^{(s-1)^{k-1-j}}}\, \, , \qquad \text{for} \,\, k \, \text{\rm odd}.
\end{equation}
\end{theorem}

\noindent{\it Proof.}
We make the proof by mathematical induction over $k$.

For $k=0$, (\ref{GenSylvE}) is trivially satisfied since, in the right hand side, we have a ratio of empty products.

For $k=1$, since the product in the denominator of the right hand side is an empty product, we have
$$
\det\, M_{1}  \cdot  \left[ \det\, M_{0} \right]^{s-1}=   \det B_{0}, 
$$
where $\displaystyle M_{1} = M \left( \!\!\! \begin{array}{c}
1 ,\ldots ,t+s \\1 , \ldots,t+s \end{array} \!\!\! \right)$, $\det M_{0} = a_{t,t}^{(t-1)}$, and  $B_0$
is the matrix with elements $b_{ij}^{(0)}= M \left[ \!\!  \begin{array}{c}
1, \ldots,  t, t+ i  \\ 1, \ldots,  t, t+ j
 \end{array}  \!\! \right]=a_{t+i,t+j}^{(t)}$, for $ 1 \le i,j\le s$. So it is exactly the
Sylvester's identity (\ref{IdSylv}) applied to $M_{1}$, and  (\ref{GenSylvO}) holds.

In the inductive step we shall prove that, if (\ref{GenSylvO}) holds for $k$ odd, then
(\ref{GenSylvE}) holds for $k+1$ even, moreover
if (\ref{GenSylvE}) holds for $k$ even, then (\ref{GenSylvO}) holds for $k+1$ odd.

Assume that (\ref{GenSylvO}) is satisfied for $k<q$ odd.
We apply (\ref{IdSylv}) to the matrix $M_{k+1}$, of order $(m_k+s)$, with
$\det M_{k}= a_{m_k,m_k}^{(m_k-1)}$ fixed, and we have
\begin{equation} \label{SyO} \det\, M_{k+1} \cdot
\left[ \det\, M_{k} \right]^{s-1}= \det B_{k} \end{equation}
where $B_k$ is the matrix with elements as in (\ref{BIJ}).
Substituting in (\ref{SyO}) the expression for $\det M_k$ obtained from
(\ref{GenSylvO}), leads to
$$
\frac{\det\, M_{k+1}} {\left[ \det\, M_{0} \right]
^{(s-1)^{k+1}}} \cdot \left[ \frac{\displaystyle \prod_{\substack{ i=0 \\ i  \text{\rm ~even}}}^{k-1}
\, \left[  \det B_{i} \right]^{(s-1)^{k-1-i}} } {\displaystyle \prod_{\substack{ j=1 \\ j   \text{\rm ~odd}}}^{k-2}
 \,\left[  \det B_{j} \right]^{(s-1)^{k-1-j}} }
\right]^{s-1}=  \det B_{k} \, ,
$$
and then we get
$$
\frac{\det\, M_{k+1}} {\left[ \det\, M_{0} \right]
^{(s-1)^{k+1}}}= \frac{\displaystyle \prod_{\substack{ j=1 \\ j
\text{\rm ~odd}}}^{k}
 \,\left[  \det B_{j} \right]^{(s-1)^{k-j}}}
{\displaystyle \prod_{\substack{ i=0 \\ i   \text{\rm ~even}}}^{k-1}
 \,\left[  \det B_{i} \right]^{(s-1)^{k-i}}}\, \,, $$
which is (\ref{GenSylvE}) for $k+1$ even.

Assume, now, that (\ref{GenSylvE}) holds for $k<q$ even.
\noindent  By using (\ref{IdSylv}), the relation (\ref{SyO}) holds again.
Substituting in (\ref{SyO}) the expression for $\det M_k$, obtained from
(\ref{GenSylvE}),  gives now
$$
\det\, M_{k+1} \cdot {\left[ \det\, M_{0} \right] ^{(s-1)^{k+1}}}
\cdot \left[ \frac{\displaystyle \prod_{\substack{ i=1 \\ i \text{\rm ~odd}}}^{k-1} \, \left[
 \det B_{i} \right]^{(s-1)^{k-1-i}} } { \displaystyle \prod_{\substack{ j=0 \\ j  \text{\rm ~even}}}^{k-2}
 \,\left[  \det B_{j} \right]^{(s-1)^{k-1-j}} }
\right]^{s-1}=  \det B_{k},
$$
and then
$$
\det\, M_{k+1} \cdot {\left[ \det\, M_{0} \right]
^{(s-1)^{k+1}}}= \frac{\displaystyle \prod_{\substack{ j=0 \\ j
\text{\rm ~even}}}^{k}
 \,\left[  \det B_{j} \right]^{(s-1)^{k-j}}}
{\displaystyle \prod_{\substack{ i=1 \\ i   \text{\rm ~odd}}}^{k-1}
 \,\left[  \det B_{i} \right]^{(s-1)^{k-i}}} \, ,
$$
which proves (\ref{GenSylvO}) for $k+1$ odd. \\$\Box $

So, now, we have established a relation between determinants obtained by bordering a leading principal submatrix of order
$t$ of a matrix $M$, with blocks of $s$ rows and $s$ columns.
 $m=n$ (i.e. $s=n-t$), (\ref{GenSylvO}) becomes
 exactly the classical Sylvester's determinantal identity (\ref{IdSylv}).

\noindent{\it Remark:} When $s=1$, we have $\det M_k=\det B_{k-1}= a_{t+k,t+k}^{(t+k-1)}$, for $ 0 \leq k \leq \min(n-t,m-t)$,
and the $B_{k-1}$ are scalars.

\subsection{Applications}

Obviously, from the computational point of
view, if $s$ is large, the formulae given in the Theorem \ref{th}
are not very interesting.
So, let $M$ a square  matrix of order $n$, and consider the particular case where $s=2$, chosen $t$ such that
$n-t$ is even, $q = (n-t)/2$, even or odd.
We set
$$ M =
\left( \begin{array}{cc}
A_{11} & A_{12} \\
A_{21}^T & A_{22}
\end{array}
\right), $$
where  $A_{11}$ is a  square matrix of order $t$, $A_{22}$ is a square matrix
of order $(n-t)=2q$, $A_{12}, A_{21}$ matrices of dimension $t \times (n-t)= t \times 2q$.

Formulae (\ref{GenSylvE}) and (\ref{GenSylvO}), for $k=q$, become
\begin{equation} \label{GenSylvEA} \frac{\det\, M} {  \det
A_{11}  }= {\displaystyle
\prod_{\substack{ i=1 \\ i  \, \text{odd}}}^{q-1}
  \det B_{i}} \;\big/
{\displaystyle \prod_{\substack{ j=0 \\ j   \text{~even}}}^{q-2}
  \det B_{j} }\,\,, \qquad  \,\,\,\,\,\,\,\,\,\,\,\,\,\,\,\,\,\,\,\,\,\, q \text{~even} \, ,
\end{equation}
\begin{equation} \label{GenSylvOA} \det\, M \, \cdot   \det
A_{11} = {\displaystyle \prod_{\substack{
i=0 \\ i   \text{~even}}}^{q-1}
   \det B_{i}} \; \big/
{\displaystyle \prod_{\substack{ j=1 \\ j   \text{~odd}}}^{q-2}
   \det B_{j} }\, \, , \qquad q  \text{~odd},
\end{equation}
where the $B_i$ matrices, computed as in (\ref{BIJ}), are all of order $2$.

Consider the particular case where $q=1$. Let $D=A_{11}$ be a
square matrix of order $n-2$,
$$ A_{22}= \left(
\begin{array}{cc}
a & e \\
g & l
\end{array}
\right)  , \qquad A_{12}= \left(
\begin{array}{cc}
{\mathbf{c}} & {\mathbf{f}}
\end{array}
\right), \qquad
A_{21}= \left(
\begin{array}{cc}
{\mathbf{b}} & {\mathbf{h}}
\end{array}
\right),
$$
with $a$, $e$, $g$, $l$ scalars, and ${\mathbf{b}},
{\mathbf{c}}, {\mathbf{f}}, {\mathbf{h}}$   column
vectors of  length $(n-2)$.
Formula (\ref{GenSylvOA}) gives
$$
\det  M_1  \cdot \det  D = \det B_{0},
$$
and since, in this case,
$$\det B_0=
\left| \begin{array}{cccc}
\left| \begin{array}{cc} D & {\mathbf{c}}
 \\ {\mathbf{b}}^{T} & a
\end{array}
\right| &
\left| \begin{array}{cc}
D & {\mathbf{f}} \\ {\mathbf{b}}^{T} & e
\end{array}
\right| \\ \\
 \left| \begin{array}{cc}
D & {\mathbf{c}} \\ {\mathbf{h}}^{T} & g
\end{array}
\right| &
\left| \begin{array}{cc}
D & {\mathbf{f}} \\ {\mathbf{h}}^{T} & l
\end{array}
\right|
\end{array}
\right|,
$$
we recover exactly the rule (\ref{SQUARE}).

Let now $M^\prime$ a square  matrix of order $n$, partitioned as
\begin{equation}
\label{Mprimo} M^{\prime} =
\left( \begin{array}{ccc}
A & B^T & E \\
C & D & F\\
G^T & H^T & L
\end{array}
\right),
\end{equation}
where the  block $D$ is a square  matrix of order $t$,
$A,E,G,L$ are square matrices
of order $2$, and $B,C,F,H$ are of dimension $t \times 2$, and $n=t+4$.
Setting
$$ M= \left( \begin{array}{ccc}
D & C & F \\
B^T & A & E\\
H^T & G^T & L
\end{array}
\right), $$
we have $\det M^{\prime} = \det M$. So considering $s=q=2$,  (\ref{GenSylvEA}) gives
\begin{equation}
\label{GenSylvEAnew}
 \frac{\det M } {\det D}=\frac{ \det B_{ 1}}
{ \det B_{ 0} },
\end{equation}
with $B_0$ and $B_1$ of order $2$.

For the above matrices we denote
\begin{eqnarray*}
&& B=  \left(
\begin{array}{cc}
{\mathbf{b}_1} & {\mathbf{b}_2}
\end{array}
\right), \;
C=  \left(
\begin{array}{cc}
{\mathbf{c}_1} & {\mathbf{c}_2}
\end{array}
\right), \;
E=  \left(
\begin{array}{cc}
{\mathbf{e}_1} & {\mathbf{e}_2}
\end{array}
\right), \;
F=  \left(
\begin{array}{cc}
{\mathbf{f}_1} & {\mathbf{f}_2}
\end{array}
\right)\!, \;\\
&&
G=  \left(
\begin{array}{cc}
{\mathbf{g}_1} & {\mathbf{g}_2}
\end{array}
\right)\!, \; H=  \left(
\begin{array}{cc}
{\mathbf{h}_1} & {\mathbf{h}_2}
\end{array}
\right), \;  
 \\
&&A= (a_{ij}), \qquad  L= (l_{ij}),\; \mbox{\rm ~~~~~~~for~}i,j = 1,2,
\end{eqnarray*}
where $\mathbf{b}_i,\mathbf{c}_i,\mathbf{e}_i,\mathbf{f}_i,\mathbf{g}_i,\mathbf{h}_i,  i = 1,2$, are column vectors.
We have
$$\det B_0=
\left| \begin{array}{cccc}
\left| \begin{array}{cc} D & {\mathbf{c}_1}
 \\ {\mathbf{b}_1}^{T} & a_{11}
\end{array}
\right| &
\left| \begin{array}{cc} D & {\mathbf{c}_2}
 \\ {\mathbf{b}_1}^{T} & a_{12}
\end{array}
\right|\\ \\
\left| \begin{array}{cc} D & {\mathbf{c}_1}
 \\ {\mathbf{b}_2}^{T} & a_{21}
\end{array}
\right| &
\left| \begin{array}{cc} D & {\mathbf{c}_2}
 \\ {\mathbf{b}_2}^{T} & a_{22}
\end{array}
\right|
\end{array}
\right|=
\left| \begin{array}{cccc}
\left| \begin{array}{cc} a_{11} & {\mathbf{b}_1}^T
 \\ {\mathbf{c}_1} & D
\end{array}
\right| &
\left| \begin{array}{cc} a_{12} & {\mathbf{b}_1}^T
 \\ {\mathbf{c}_2} & D
\end{array}
\right|\\ \\
\left| \begin{array}{cc} a_{21} & {\mathbf{b}_2}^T
 \\ {\mathbf{c}_1} & D
\end{array}
\right| &
\left| \begin{array}{cc} a_{22} & {\mathbf{b}_2}^T
 \\ {\mathbf{c}_2} & D
\end{array}
\right|
\end{array}
\right|,
$$
$$
\det B_1=
\left| \begin{array}{cc}
\left|
  \begin{array}{@{~}c@{~}c@{~}c@{~}} D & C &{\mathbf{f}_1}
 \\ B^T & A & {\mathbf{e}_1}\\
 {\mathbf{h}_1}^T
              & {\mathbf{g}_1^T}
 & l_{11}
\end{array}
\right| &
\left|  \begin{array}{@{~}c@{~}c@{~}c@{~}} D & C & {\mathbf{f}_2}
 \\ B^T & A & {\mathbf{e}_2}\\
{\mathbf{h}_1}^T
              & {\mathbf{g}_1^T}
 & l_{12}
\end{array}
\right|
 \\ \\
\left| \begin{array}{@{~}c@{~}c@{~}c@{~}} D & C & {\mathbf{f}_1}
 \\ B^T & A & {\mathbf{e}_1}\\
 {\mathbf{h}_2}^T
              & {\mathbf{g}_2^T}
 & l_{21}
\end{array}
\right| &
\left|  \begin{array}{@{~}c@{~}c@{~}c@{~}} D & C & {\mathbf{f}_2}
 \\ B^T & A & {\mathbf{e}_2} \\
{\mathbf{h}_2}^T
              & {\mathbf{g}_2^T}
 & l_{22}
\end{array}
\right|
\end{array}
\right|=
\left| \begin{array}{cc}
\left|
  \begin{array}{@{~}c@{~}c@{~}c@{~}} A & B^T &{\mathbf{e}_1}
 \\ C & D & {\mathbf{f}_1}\\
 {\mathbf{g}_1^T}
              & {\mathbf{h}_1}^T
 & l_{11}
\end{array}
\right| &
\left|  \begin{array}{@{~}c@{~}c@{~}c@{~}} A & B^T & {\mathbf{e}_2}
 \\ C & D &  {\mathbf{f}_2}\\
{\mathbf{g}_1^T}
              & {\mathbf{h}_1}^T
 & l_{12}
\end{array}
\right|
 \\ \\
\left| \begin{array}{@{~}c@{~}c@{~}c@{~}} A & B^T  & {\mathbf{e}_1}
 \\ C & D & {\mathbf{f}_1}\\
 {\mathbf{g}_2^T}
              & {\mathbf{h}_2}^T
 & l_{21}
\end{array}
\right| &
\left|  \begin{array}{@{~}c@{~}c@{~}c@{~}} A & B^T  & {\mathbf{e}_2}
 \\ C & D & {\mathbf{f}_2} \\
{\mathbf{g}_2^T}
              & {\mathbf{h}_2}^T
 & l_{22}
\end{array}
\right|
\end{array}
\right|.
$$
So this application is  a generalization of the {\it Case 2} of the Section \ref{SEC2}.
\subsubsection{Examples}
\noindent{\bf Example 1:}
Let $n=10$, $t=s=2$, $q=(n-t)/{2}=4$  even. Thanks to
(\ref{GenSylvEA}),
we have
$$
\frac{\det  M} {\det\, A_{11}}= \frac{ \det B_{1} \cdot  \det B_{
3}} { \det B_{ 0} \cdot  \det B_{ 2}}.
$$
If $n=10$, $t=6$, $s=2$, $q=(n-t)/{2}=3$  odd, thanks to
(\ref{GenSylvOA}),
we have
$$
\det  M \cdot \det A_{11}= \frac{ \det B_{0} \cdot  \det B_{
2}} { \det B_{ 1} }.
$$

\noindent{\bf Example 2:}
We consider a square  matrix $M^\prime$ of order $n=8$ partitioned like in (\ref{Mprimo}), with
$\displaystyle D=
M^\prime \! \left( \!\!\! \begin{array}{l}
3, \ldots, 6\\
3, \ldots, 6
\end{array}  \!\!\!\right)$. Then $t=4, s=q=2$, and
\begin{eqnarray*}
&&\displaystyle A=
M^\prime \! \left( \!\!\! \begin{array}{l}
1,2\\
1,2
\end{array}  \!\!\!\right), \qquad E=
M^\prime \! \left( \!\!\! \begin{array}{l}
1,2\\
7,8
\end{array}  \!\!\!\right), \qquad
G=
M^\prime \! \left( \!\!\! \begin{array}{l}
7,8\\
1,2
\end{array}  \!\!\!\right), \qquad
L=
M^\prime \! \left( \!\!\! \begin{array}{l}
7,8\\
7,8
\end{array}  \!\!\!\right),\\
&&
B^T=
M^\prime \! \left( \!\!\! \begin{array}{c}
1,2\\
3, \ldots, 6
\end{array}  \!\!\!\right)\!,\,
H^T=
M^\prime \! \left( \!\!\! \begin{array}{c}
7,8\\
3, \ldots, 6
\end{array}  \!\!\!\right)\!,\,
C=
M^\prime \! \left( \!\!\! \begin{array}{c}
3, \ldots, 6\\
1,2
\end{array}  \!\!\!\right)\!,\,
F=
M^\prime \! \left( \!\!\! \begin{array}{c}
3, \ldots, 6\\
7,8
\end{array}  \!\!\!\right).
\end{eqnarray*}
So we can apply (\ref{GenSylvEAnew}), and we obtain
$$
 \frac{\det\, M^\prime  } {\det\, D}=
\frac{ M^\prime  \left[ \!\!\! \begin{array}{c}
1 ,\ldots ,6,7 \\1 ,\ldots, 6,7
 \end{array} \!\!\! \right]  \cdot
 M^\prime   \left[ \!\!\! \begin{array}{c}
1,\ldots, 6, 8 \\1,\ldots, 6, 8
 \end{array} \!\!\! \right]-
 M^\prime \left[ \!\!\! \begin{array}{c}
1,\ldots, 6, 7 \\1,\ldots, 6, 8
 \end{array} \!\!\! \right]  \cdot
 M^\prime   \left[ \!\!\! \begin{array}{c}
1,\ldots, 6, 8 \\1,\ldots, 6, 7
 \end{array} \!\!\! \right]}
 {M^\prime  \left[ \!\!\! \begin{array}{c}
1,3,\ldots, 6 \\1,3,\ldots, 6
 \end{array} \!\!\! \right]  \cdot
 M^\prime   \left[ \!\!\! \begin{array}{c}
2,3,\ldots, 6 \\2,3,\ldots, 6
 \end{array} \!\!\! \right]-
 M^\prime  \left[ \!\!\! \begin{array}{c}
1,3,\ldots, 6 \\2,3,\ldots, 6
 \end{array} \!\!\! \right] \cdot
 M^\prime   \left[ \!\!\! \begin{array}{c}
2,3,\ldots, 6 \\1,3,\ldots, 6
 \end{array} \!\!\! \right]} .
$$

\section{Conclusions}
In this paper we presented several generalizations of the
 Sylvester's determinantal identity, proposed from various authors, and described here in  a unified way.
  We started to study those formulae because the classical formula is a very important
  tool for finding extrapolation algorithms \cite{THEBOOK,schur} for accelerating the convergence of scalar
  and vector sequences, and
  we need to have at our disposal generalizations able to be used in building new extrapolation algorithms for vector sequences
  (a case more difficult) \cite{BR1}. During our study, we had the idea for the new generalization proposed. 

\section*{Remark} This is an updated version of a work presented in two congresses in 2008.


\begin{thebibliography}{99}

\bibitem{Aitk}
A.C. Aitken
\newblock {\em Determinants and Matrices},
\newblock Oliver and Boyd, Edinburgh, 1964.


\bibitem{AAM}
A.G. Akritas, E.K. Akritas, G.I. Malaschonok,
\newblock {\em Various proofs of Sylvester's (determinant) identity},
\newblock Math. Comput. Simulation 42 (1996)  585--593.



\bibitem{Bair}
E.H. Bareiss,
\newblock {\em Sylvester's identity and multistep integer-preserving Gaussian elimination},
\newblock Math. Comp. 22 (1968) 565--578.

\bibitem{Bair2}
E.H. Bareiss,
\newblock {\em Computational solution of matrix problems over an integral domain},
\newblock J. Inst. Math. Appl. 10 (1972) 68--104.

\bibitem{BeckLaba}
B. Beckermann, G. Labahn,
\newblock {\em Fraction-free computation of matrix rational interpolants
and matrix GCD's},
\newblock SIAM J. Matrix Anal. Appl. 22 (2000) 114--144.

\bibitem{BeckMuhl}
B. Beckermann, G. M\"uhlbach,
\newblock {\em A general determinantal identity of Sylvester's type and some applications},
\newblock Linear Algebra Appl. 197 (1994) 93--112.


\bibitem{BR}
C. Brezinski,
\newblock {\em A general extrapolation algorithm},
\newblock Numer. Math.  35  (1980) 175--187.



\bibitem{BR1}
C. Brezinski,
\newblock {\em  Some determinantal identities in a vector space, with applications},
\newblock in: H. Werner,
H.J. B\"{u}nger (Eds.), Pad\'e Approximation and its Applications, Bad Honnef 1983, Lecture Notes in
Mathematics, vol. 1071, Springer-Verlag, Berlin, 1984, pp. 1–-11.

\bibitem{schur}
C. Brezinski,
\newblock {\em  Schur complements and applications in numerical analysis},
in: F. Zang (Ed.), The Schur complement and its applications, Springer, New York, 2005, pp. 227--258.

\bibitem{THEBOOK}
C. Brezinski, M. Redivo-Zaglia,
\newblock {\em  Extrapolation Methods. Theory and Practice},
\newblock Studies in Computational Mathematics, North-Holland, Amsterdam, 1991.



\bibitem{Br-Rz}
C. Brezinski, M. Redivo-Zaglia,
\newblock {\em A Schur complement approach to a general
extrapolation algorithm},
\newblock Linear Algebra Appl. 368 (2003)  279--301.


\bibitem{BrSch}
R.A. Brualdi, H. Schneider,
\newblock {\em Determinantal identities: Gauss, Schur, Cauchy, Sylvester, Kronecker, Jacobi, Binet, Laplace, Muir, and Cayley.},
\newblock Linear Algebra Appl. 52/53 (1983)  769--791.



\bibitem{chio}
H. Eves,
\newblock {\em Chi\'o Expansion},
\newblock  in: Elementary Matrix Theory.  Dover, New York, 1996, pp. 129--136.




\bibitem{gant}
F.R. Gantmacher,
\newblock {\em The Theory of Matrices},
\newblock   Chelsea Publishing Co., New York, 1959.


\bibitem{GaLoRam}
 M. Gasca, A. Lopez-Carmona, V. Ramirez,
\newblock {\em A generalized Sylvester's identity on determinants and its applications to interpolation problems},
\newblock  in: Schempp, Zeller (Eds), Multivariate Approximation Theory II,  Birkhauser V., 1982, pp. 171--184.

\bibitem{Konva}
M. Konvalinka,
\newblock {\em Non-commutative Sylvester's determinantal identity},
\newblock Electron. J. Combin.  14 (2007), R42, 29 pp. (electronic).

\bibitem{Leclerc}
B. Leclerc,
\newblock {\em On identities satisfied by minors of a matrix},
\newblock  Adv.  Math. 100 (1993)  101--132.

\bibitem{Muir}
T. Muir,
\newblock {\em The law of extensible minors in determinants},
\newblock  Trans. Roy. Soc. Edinburgh 30 (1881) 1--4.


\bibitem{Muhl}
G. M\"uhlbach,
\newblock {\em On extending determinantal identities},
\newblock  Linear Algebra Appl. 132 (1990)  145--162.


\bibitem{MuhlGas}
G. M\"uhlbach, M. Gasca,
\newblock {\em A generalization of Sylvester's identity on determinants and some applications},
\newblock  Linear Algebra Appl.  66 (1985)  221--234.



\bibitem{Mulders}
T. Mulders,
\newblock {\em A generalized Sylvester identity and fraction-free random Gaussian elimination},
\newblock  Journal Symbolic Comput.    31 (2001) 447--460.

\bibitem{Sylv}
J.J. Sylvester,
\newblock {\em On the relation between the minor determinants of linearly equivalent quadratic functions},
\newblock  Philos. Mag., (4th. Series)  1 (1851) 295--305.

\bibitem{Yak1}
O.S. Yakovlev,
\newblock {\em Some determinant identities},
\newblock  Ukrainian Math. J.     26 (1974) 352--355.


\bibitem{Yak2}
O.S. Yakovlev,
\newblock {\em A generalization of determinant identities},
\newblock  Ukrainian Math. J. 30 (1978) 643--647.


\end{thebibliography}
\end{document}